\documentclass{article}

\usepackage{amsfonts}
\usepackage{amsmath}
\usepackage{amsthm}
\usepackage{mathtools}
\usepackage{cite}
\usepackage{yfonts}
\usepackage[linesnumbered,boxed]{algorithm2e}
\usepackage{tikz-cd}
\usepackage{hyperref}
\usepackage{url}
\usepackage{cleveref}

\Crefname{ALC@unique}{Line}{Lines}
\newcounter{myalg}
\AtBeginEnvironment{algorithmic}{\refstepcounter{myalg}}
\makeatletter
\@addtoreset{ALC@unique}{myalg}
\makeatother

\newtheorem{theorem}{Theorem}[section]
\newtheorem{lemma}[theorem]{Lemma}
\newtheorem{prop}[theorem]{Proposition}
\newtheorem{cor}[theorem]{Corollary}
\newtheorem{fact}[theorem]{Fact}
\newtheorem{definition}[theorem]{Definition}

\theoremstyle{remark}
\newtheorem{remark}[theorem]{Remark}

\DeclareMathOperator{\Tr}{Tr}
\DeclareMathOperator{\Diag}{Diag}
\DeclareMathOperator{\Gal}{Gal}

\DeclareMathOperator{\Ker}{Ker}
\DeclareMathOperator{\im}{Im}
\DeclareMathOperator{\Cl}{Cl}

\DeclareMathOperator{\Dc}{\mathcal{D}}
\DeclareMathOperator{\Pc}{\mathcal{P}}
\DeclareMathOperator{\cl}{Cl}
\DeclareMathOperator{\Pl}{Pl}
\DeclareMathOperator{\Br}{Br}
\DeclareMathOperator{\ord}{ord}
\DeclareMathOperator{\Supp}{Supp}

\newcommand{\Q}{\mathbb{Q}}
\newcommand{\Z}{\mathbb{Z}}
\newcommand{\N}{\mathbb{N}}
\newcommand{\F}{\mathbb{F}}

\newcommand{\Oc}{\mathcal{O}}

\newcommand{\p}{\mathfrak{p}}
\newcommand{\1}{\mathbf{1}}
\newcommand{\mpp}{\mathfrak{m}}
\renewcommand{\epsilon}{\varepsilon}

\title{Efficient computations in central simple algebras using Amitsur cohomology}
\author{Péter Kutas\thanks{Faculty of Informatics, Eötvös Loránd University and School of Computer Science, University of Birmingham (\url{https://sites.google.com/view/peterkutas89}). The first author is supported by the Hungarian Ministry of Innovation and Technology NRDI Office within the framework of the Quantum Information National Laboratory Program, the János Bolyai Research Scholarship of the Hungarian Academy of Sciences and by the UNKP-22-5 New National Excellence Program. The first author is also partly supported by EPSRC through grant number EP/V011324/1.} \and Mickaël Montessinos\thanks{Institute of Mathematics, Faculty of Mathematics and Informatics, Vilnius University (\url{http://mickael.montessinos.fr}).}}
\begin{document}

\maketitle

\begin{abstract}
    We introduce a presentation for central simple algebras over a field \(k\) using Amitsur cohomology. We provide efficient algorithms for computing a cocycle corresponding to any such algebra given by structure constants. If \(k\) is a number field, we use this presentation to prove that the explicit isomorphism problem (i.e., finding an isomorphism between central simple algebras given by structure constants) reduces to \(S\)-unit group computation and other related number theoretical computational problems. This also yields, conditionally on the generalised Riemann hypothesis, the first polynomial quantum algorithm for the explicit isomorphism problem over number fields.
\end{abstract}

\section{Introduction}
The \emph{explicit isomorphism problem} is the algorithmic problem of, given a field \(k\) and a \(k\)-algebra \(A\) isomorphic to \(M_d(k)\), constructing an explicit isomorphism \(\varphi\colon A \to M_d(k)\). The explicit isomorphism problem may be thought of as a natural problem in computational representation theory. Given a \(k\)-algebra \(A\), one may wish to assay it. That is, compute the Jacobson radical of \(A\), and the decomposition of the semi-simple quotient of \(A\) as a sum of simple \(k\)-algebras, themselves identified to some \(M_d(D)\), for \(D\) a division \(k\)-algebra. In general, the hard part of this task is to find an isomorphism \(A \to M_n(D)\) when \(A\) is simple. A general recipe for solving this problem is to identify the Brauer class of \(D\) over its centre \(K/F\), find structure constants for \(M_n(D^{op})\) and then compute an explicit isomorphism \(A \otimes M_n(D^{op}) \simeq M_m(K)\) \cite{ivanyos2012splitting,gomez2022primitive,csahok2022explicit}.

Applications of the explicit isomorphism problem go beyond the mere computational theory of finite-dimensional associative algebras. In arithmetic geometry, the problem is relevant for trivialising obstruction algebras in explicit descent over elliptic curves \cite{cremona2008explicit,cremona2009explicit,cremona2015explicit,fisher2013explicit} and computation of Cassel-Tate pairings \cite{fisher2014computing,yan2021computing}. The problem is also connected to the parametrisation of Severi-Brauer surfaces \cite{de2006lie}. Recent work in algebraic complexity theory reduced the determinant equivalence test to the explicit isomorphism problem \cite{garg2019determinant}. Finally, the explicit isomorphism problem over the function field \(\F_q(T)\) is also relevant to error correcting codes \cite{gomez2022primitive}.

It is well known that the Brauer group of equivalence classes of central simple \(k\)-algebras may be described as the second Galois cohomology group \(H^2(k,k_{sep}^\times)\). The remarkable fact that these two groups are isomorphic may be seen as a consequence of the crossed-product presentation of central simple algebras containing a maximal commutative subalgebra which is a Galois extension of the base field.

When a cohomological presentation is known for a certain algebra, the explicit isomorphism problem translates into a multiplicative algebraic equation. Indeed, the algebra is represented by a cocycle, and since the algebra is split, the cocycle is a coboundary. Then, an explicit isomorphism to a matrix algebra may be computed using a cochain whose differential is this coboundary. When the base field is a number field and a presentation of the algebra (either a cyclic or a crossed product presentation) is known, the resulting equation may be solved by computing a certain group of \(S\)-units \cite{simon2002solving,fieker2009minimizing}. 

Before Noether introduced crossed-product algebras, Brauer had already introduced a factor set presentation which allowed to describe an algebra with the knowledge of any maximal commutative subalgebra, not necessarily a Galois extension of the base field. We refer the reader to \cite{jacobson1996finite} for a modern description of presentations introduced in the works of Brauer and Noether. It is interesting to note that Adamson developed a cohomology theory for non-Galois field extensions which describes Brauer factor sets as cocycles \cite{adamson1954cohomology}, much like classical Galois cohomology describes Noether factor sets as cocycles.

A line of work introduced Azumaya algebras, which generalise the notion of central simple algebra (and therefore of the Brauer group) over arbitrary rings \cite{azumaya1951maximally,auslander1960brauer}, and then later over schemes (see \cite{colliot2021brauer} for a modern treatment). While trying to connect the usual Galois cohomological description of the Brauer group and another description introduced by Hochschild, Amitsur came upon an exact complex which yields a cohomology theory suitable to classify central simple algebras, the so-called Amitsur cohomology. As the complex he defined may be constructed over arbitrary rings, a line of work generalised his classification result to Azumaya algebras over increasingly general classes of ring extensions \cite{amitsur1959simple,rosenberg1960amitsur,chase1965amitsur}.

Cohomological presentations have the potential of being used both for a computational representation of a central simple algebra and for solving the explicit isomorphism problem. For instance, in the computer algebra system PARI/GP \cite{pari}, central simple algebras over a number field may be represented as cyclic algebras. While the already existing constructions relying on Amitsur cohomology are abstract and do not lend themselves easily to a practical implementation, the crossed product and cyclic presentations have the inconvenient that one is required to know a Galois, or even cyclic, maximal subfield in order to convert a given algebra into these presentations. While such a subfield always exists for algebras over global fields (as per the Brauer-Hasse-Noether theorem), there is no known efficient algorithm to find one (except for quaternion and degree 3 central simple algebras) to the best of our knowledge. The presentation using Brauer factor sets may be constructed using an arbitrary maximal commutative subalgebra, which is easy to find, but factor sets take value in a normal splitting field of the subalgebra used. Results in arithmetic statistics \cite{eberhard2022characteristic} suggest that a random element of such an algebra will generate a maximal subfield whose Galois group is the full symmetric group, and hence will have a Galois closure of degree \(n!\) (where \(n\) is the degree of the algebra), preventing it to be used efficiently in a computation.

\subsection{Our contributions}
In \Cref{Sec:Amitsur}, we give a presentation for central simple algebras over a field using Amitsur cocycles. Our presentation only requires knowledge of a separable maximal commutative subalgebra, which may be found efficiently in deterministic polynomial time \cite{graaf2000finding}. Furthermore, multiplication in our representation may be computed using \(O(n^3)\) base field multiplications, making it as efficient as the naive representation using structure constants. In \Cref{Sec:AmitsurAlgebra} we construct a central simple algebra from an Amitsur \(2\)-cocycle. In \Cref{Sec:CompRepAmi}, we discuss the computational representation of Amitsur cocycles and related algebras. Finally, in \Cref{Sec:AmitsurAlgebraConstruction} we give an efficient algorithm for computing a cocycle representing a given central simple algebra.

In \Cref{Sec:Trivialisation}, we prove that over a number field, an Amitsur \(1\)-cochain whose differential is a given \(2\)-coboundary may be found within a certain group of \(S\)-units. We then prove that, assuming the generalised Riemann hypothesis (henceforth denoted by \emph{GRH}), the explicit isomorphism problem is solved by a polynomial quantum algorithm.
\subsection{Related works}
Below, we review existing classical algorithms for solving the explicit isomorphism problem over various fields. 

In the case of a finite base field, a polynomial-time algorithm was introduced by Ronyai in \cite{ronyai1990computing}.

Instances of the problem for \(\Q\)-algebras isomorphic to \(M_n(\Q)\) were first treated separately for small values of \(n\). When \(n=2\), the problem reduces to finding a rational point on a projective conic \cite[Theorem 5.5.4]{voight2021quaternion}, which is solved for instance in \cite{cremona2003efficient}. Then, \cite{de2006lie} presented a subexponential algorithm when \(n=3\) by finding a cyclic presentation and solving a cubic norm equation. The case \(n=4\) is tackled in \cite{pilnikova2007trivializing} by reducing to the case of quaternion algebras over \(Q\) and \(\Q(\sqrt{d})\) and solving a norm equation.

In \cite{cremona2015explicit} an algorithm was given and studied mostly for the cases \(n=3\) and \(n=5\). It was then generalised in \cite{ivanyos2012splitting,ivanyos2013improved} to a \(K\)-algebra isomorphic to \(M_n(K)\), where \(n\) is a natural number and \(K\) is a number field. The complexity of this last algorithm is polynomial in the size of the structure constants of the input algebra, but depends exponentially on \(n\), the degree of \(K\) and the size of the discriminant of \(K\).

In 2018, \cite{ivanyos2018computing} exhibited a polynomial-time algorithm for algebras isomorphic to \(M_n(\F_q(T))\).

For the case of fixed \(n\) and varying base field, \cite{fisher2017higher,kutas2019splitting} independently gave an algorithm for an algebra isomorphic to \(M_2(\Q(\sqrt{d}))\) with complexity polynomial in \(\log(d)\). A similar algorithm for quadratic extensions of \(\F_q(T)\) was given in \cite{ivanyos2019explicit} for the case of odd \(q\). The case of even \(q\) was then treated in \cite{kutas2022finding}.

\section{Preliminaries}\label{Sec:Prelim}
\subsection{Algorithmic preliminaries}\label{Sec:AlgoPrelim}

    \subsubsection{Computational model}
        The theoretical results of \Cref{Sec:Amitsur} are valid for any base field \(k\). In order to translate them into efficient algorithms, it is necessary to have an encoding of the elements of \(k\), for which the usual field operations and linear algebraic tasks may be performed efficiently. Following \cite{graaf2000finding}, we say that such a field admits \emph{efficient linear algebra}. It is well known that finite fields and global fields, among others, admit efficient linear algebra.
        
    \subsubsection{Computational representations of algebras}\label{Sec:Presentations}
        A \(k\)-algebra may be given as input to a program in the form of a table of structure constants. Let \(A\) be a finite dimensional \(k\)-algebra of dimension \(n\), let \(V\) be the underlying vector space of \(A\) and let \((e_1,\hdots,e_n)\) be a basis of \(V\). Multiplication in \(A\) is a bilinear map which may be represented as a tensor \(\lambda \in V^\vee \otimes V^\vee \otimes V\). The table of structure constants of \(A\) for the basis \((e_1,\hdots,e_n)\) is then the coordinates of \(\lambda\) in the basis \(\left(e_i^\vee \otimes e_j^\vee \otimes e_k\right)\) of \(V^\vee \otimes V^\vee \otimes V\), where \(\left(e_i^\vee\right)\) is the basis of \(V^\vee\) dual to \((e_i)\).

        An \(F\)-algebra \(A \cong M_d(k)\) may be presented as a cyclic algebra \cite[Section 2.5]{gille2017central} or as a crossed-product algebra \cite[Section 2.6]{jacobson1996finite}. It follows readily from the algebraic constructions that, computationally, knowing a cyclic (resp. crossed-product) presentation of a central simple \(F\)-algebra \(A\) of degree \(d\) is equivalent to knowing an embedding \(F \to A\), where \(F/k\) is a cyclic (resp. Galois) extension of  \(k\) of degree \(d\).

        Field extensions and étale algebras admit a presentation as a quotient ring \(k[x]/(P(x))\), for some \(P \in k[x]\). As a result, their elements are represented as residue classes of polynomials. The question of tensor products of such algebras is discussed below in \Cref{Sec:CompRepAmi}.
    
\subsubsection{Computing rings of integers and factoring ideals}
    The ability to factor integers in polynomial time yields polynomial quantum algorithms for several tasks of computational number theory. We list some auxiliary results that will be needed later on: 
    \begin{fact}\label{fact:maxord}
    There exists a polynomial-time quantum algorithm for computing the maximal order in number fields.   
    \end{fact}
    \begin{proof}
     Zassenhaus' round 2 algorithm \cite[Algorithm 6.1.8]{cohen2013course} computes the maximal order of a number field in polynomial time, provided an oracle for factoring. Thus, one may deduce a polynomial quantum algorithm for computing maximal orders in number fields.
   
    \end{proof}
    \begin{fact}\label{fact:factoring}
    There exists a polynomial-time quantum algorithm \cite[Algorithm 2.3.22]{cohen2012advanced} for factoring ideals in number fields.    
    \end{fact}

\subsubsection{A polynomial quantum algorithm for computing class groups and groups of \(S\)-units in number fields}\label{subsec:quantSunit}
The class group of $L$, denoted by $Cl(L)$ is the group of fractional ideals modulo principal ideals. The units in $\Z_L$ form a group which is called the unit group and is denoted by $U_L$. All these are important objects in number theory and they motivate the following three algorithmic problems: 
\begin{enumerate}
    \item Compute generators of $Cl(L)$ and their relations,
    \item Compute generators of $U_L$,
    \item Given a principal ideal $I$, find a generator element of $I$.
\end{enumerate}

To the best of our knowledge, there only exist polynomial-time classical algorithms for these problems in very special cases. In \cite{biasse2014subexponential} the authors propose a subexponential heuristic algorithm for computing class groups and unit groups for arbitrary degree number fields. In \cite{bauch2017short} and \cite{biasse2017computing} the authors propose subexponential heuristic algorithms for the principal ideal problem for a certain subclass of number fields (multiquadratic extensions of $\mathbb{Q}$ and power-of-two cyclotomic fields). 

Hallgren \cite{hallgren2005fast} first proposed a polynomial-time quantum algorithm for computing the class group and unit group for bounded degree number fields, where GRH is assumed for computing generators of the class group. Later, polynomial-time algorithms were proposed (under GRH for the computation of the class group) for these tasks for arbitrary number fields \cite{eisentrager2014quantum},\cite{biasse2016efficient}. Furthermore, \cite{biasse2016efficient} also provides a polynomial-time quantum algorithm for computing generators of principal ideals (dubbed the principal ideal problem). In the class group computation, GRH is required to provide a polynomially sized generating set of the class group. Bach \cite{bach1990explicit} proved that, assuming GRH, the class group can be generated by prime ideals whose norm is at most $48\log^2 |\Delta_K|$ where $\Delta_K$ is the discriminant of the number field. We summarise this discussion into the following facts:
\begin{fact}[{\cite[Fact 4.1]{biasse2016efficient}}]\label{fact:Bach}
Let $K$ be a number field and $\Delta_K$ its discriminant. Under GRH, there exists a polynomial-time algorithm computing a generator set of $Cl(K)$ of size polynomial in $\log(|\Delta_K|)$.
\end{fact}
\begin{proof}
Using the Bach bound one simply enumerates all prime ideals of norm bounded by $12\log^2(|\Delta_K|)$. The number of such prime ideals is polynomial in \(\log(|\Delta_K|)\).
\end{proof}
\begin{fact}\label{fact:class group}
Assuming GRH, there exists a polynomial-time quantum algorithm \cite{biasse2016efficient} for computing the class group of a number field.     
\end{fact}
\begin{fact}\label{fact:s-unit}[{\cite[Theorem 1.1]{biasse2016efficient}}]
    Let $S$ be a set of prime ideals of a number field $K$. There exists a polynomial-time quantum algorithm for computing generators of the group of \(S\)-units (polynomial in the size of \(S\)).  
\end{fact}

One application of $S$-unit computation (mentioned in \cite{biasse2016efficient} and \cite{simon2002solving}) is that one can compute solutions to norm equations in Galois extensions. Let $L|K$ be a Galois extension of number field and let $a\in K$. Then one would like to find a solution to  $N_{L|K}(x)=a$ (or show that it does not exist). The main result by Simon is the following: 
\begin{fact}(\cite[Theorem 4.2]{simon2002solving})\label{thm:Simon}
Let $L|K$ be a Galois extension of number fields. Let $S_0$ be a set of primes generating the relative class group $Cl_i(L|K)$. Let $S$ be a set of primes containing $S_0$ and let us denote the set of $S$-units in $L$ by $U_{L,S}$ and the set of $S$-units in $K$ by $U_{K,S}$. Then one has that: 
$$N_{L|K}(U_{L,S})=N_{L|K}(L^{*})\cap U_{K,S}$$
\end{fact}
Informally, this means that for Galois extensions one can search for $x$ in the group of $S$-units. Then by factoring $S$ and computing the class group of $L$ one can turn finding $x$ into solving a system of linear equations. 
Thus, computing $S$-units implies a polynomial quantum algorithm for solving norm equations in Galois extensions. One important issue though is that the size of $S_0$ has to polynomially bounded as the complexity could explode if $S_0$ was too large. If one assumes GRH, then \cite{bach1990explicit} implies that the size of $S_0$ is not too large. 

We summarise this discussion into the following theorem: 
\begin{theorem}\label{thm:s-unit&norm}
Let $L|K$ be a Galois extension of number fields and let $S$ be a set of primes of $L$. Then there exists a polynomial-time quantum algorithm that computes a set of generators for the group of $S$-units of $L$. Furthermore, let $a\in K$. Then there exists a polynomial-time quantum algorithm (assuming GRH) for solving the norm equation $N_{L|K}(x)=a$. 
\end{theorem}

\subsection{Brauer factor sets and central simple algebras}\label{Sec:BrauerPrelim}
For the remainder of this section, we fix a field \(k\), a separable polynomial \(P \in k[X]\) of degree \(d\) and an étale algebra \(F = k[X]/P\) over \(k\). We let \(K\) be the smallest Galois extension of \(k\) which splits \(F\), and set \(G = \Gal(K/k)\)

We will introduce definitions that are equivalent to those of Adamson cohomology \cite{adamson1954cohomology} with values in the multiplicative group. However, while Adamson cohomology is normally defined for separable field extensions, we extend our definitions to étale algebras. One difference is that cochains are indexed by \(k\)-algebra morphisms which are not necessarily injective. We may then use well-known results on Brauer factor sets to establish a cohomological description of the Brauer group (see \cite[Chapter 3]{jacobson1996finite} for a modern reference).

Let \(\Phi\) be the set of \(k\)-algebra morphisms from \(F\) to \(K\).
\begin{lemma}\label{lemma:BrauerAdamson}
     The map \(\varphi \mapsto \varphi(X)\) is a bijection between \(\Phi\) and the set of roots of \(P\) in \(K\). In particular, the set \(\Phi\) has \(d\) elements.
\end{lemma}

\begin{proof}
    If \(\varphi \in \Phi\), it is clear that \(\varphi(X)\) is a root of \(P\) in \(K\). Now, the morphism \(\varphi\) is entirely determined by \(\varphi(X)\) so the map \(\varphi \mapsto \varphi(X)\) is injective.

    Let \(r \in K\) be a root of \(P\). Then, let \(Q\) be the unique irreducible factor of \(P\) such that \(Q(r) = 0\). Then, \(F \simeq k[X]/Q \oplus k[X]/(P/Q)\) and we have a map from \(k[X]/Q \oplus k[X]/(P/Q)\) to \(k[X]/Q\) sending \((A(X),B(X))\) to \(A(r)\). This maps precomposed with the isomorphism defined above sends \(X\) to \(r\).
\end{proof}

We may then index the maps from \(F\) to \(K\), so that \(\Phi = \{\varphi_1,\hdots,\varphi_d\}\). If \(i \in [d]\) and \(g \in G\), we write \(gi\) for the index of \(g \varphi_i\).

\begin{definition}
    Let \(n \in \N\). An \(n\)-cochain is a map \[a:\begin{array}{ccl} \Phi^{n+1} &\to &K^\times \\ (\varphi_{i_0},\hdots,\varphi_{i_n}) &\mapsto &a_{i_0,\hdots,i_n}\end{array}\] which satisfies the following homogeneity condition: for \(i_0,\hdots,i_n \in [d]\), for \(g \in G\),
    \begin{equation} \label{eq:homogeneity}
        ga_{i_0,\hdots,i_n} = a_{gi_0,\hdots,gi_n}
    \end{equation}
    A cochain \(a\) may be denoted as the family \((a_{i_0,\hdots,i_n})_{i_0,\hdots,i_n \in [d]}\) of values that it takes over \(\Phi^{n+1}\).

    A cochain is said to be reduced if \(a_{i,\hdots,i} = 1\) for all \(i \in [d]\).
    Under pointwise multiplication, the \(n\)-cochains form a group which we denote by \(C^n(k,F)\). Then, the reduced cochains form a subgroup denoted by \(C'^n(k,F)\).

    If \(a \in C^n(k,F)\) and \(i_0,\hdots,i_{n+1} \in [d]\), we set \[a_{i_0,\hdots,\widehat{i_k},\hdots,i_{n+1}} \coloneqq a_{i_0,\hdots,i_{k-1},i_{k+1},\hdots,i_{n+1}}.\]

    For \(n \in \N\), there is a group homomorphism \(\Delta^n \colon C^n(k,F) \to C^{n+1}(k,F)\) defined as follows: if \(a \in C^n(k,F)\),
    \[\Delta^n(a)_{i_0,\hdots,i_{n+1}} = \prod_{k=0}^{n+1} a_{i_0,\hdots,\widehat{i_k},\hdots,i_n}^{(-1)^k}.\]

    We set \(Z^n(k,F) = \Ker(\Delta^n)\) and call elements of this subgroup \(n\)-cocycles. We write \(Z'^n(k,F)\) for the group of reduced cocycles. For \(n\) positive, we also set \(B^n(k,F) = \im(\Delta^{n-1})\) and call its elements \(n\)-coboundaries. We note that \(B^n(k,F)\) is a subgroup of \(Z^n(k,F)\). For completeness, we let \(B^0(k,F)\) be the trivial subgroup of \(Z^0(k,F)\). Finally, we set \(H^n(k,F) = Z^n(k,F)/B^n(k,F)\).

    Two cocycles \(c_1,c_2 \in Z^n(k,F)\) are called associated, denoted by \(c_1 \sim c_2\), if they have the same class in \(H^n(k,F)\). If \(b \in B^{n+1}(k,F)\), a preimage of \(b\) by \(\Delta^n\) is called a trivialisation of \(b\).
\end{definition}

\begin{remark}
    The cochain \((1)_{i,j,k \in [d]}\) is always a cocycle and is a coboundary for \(n \geq 1\). It is called the trivial cocycle.
\end{remark}

\begin{prop}\label{prop:ReducedUseless}
    For \(n\) positive, every \(n\)-cocycle is associated to a reduced cocycle.
\end{prop}

\begin{proof}
    Let \(c \in Z^n(k,F)\), and consider the cochain \(a \in C^{n-1}(k,F)\) defined by \(a_{i,\hdots,i} = c_{i,\hdots,i}^{-1}\) and \(a_{i_0,\hdots,i_{n-1}} = 1\) if the \(i_k\) are not all equal to one another. Then, the cocycle \(c \Delta^{n-1}(a)\) is reduced and associated to \(c\).
\end{proof}

We may now recall the explicit description of the isomorphism \(H^2(k,F) \simeq Br(F/k)\) (that is, the group of Brauer classes of central simple \(k\)-algebras that are split by \(F\)), as exposed in \cite[Chapter 2]{jacobson1996finite}. Using Lemma \ref{lemma:BrauerAdamson}, we observe that what is called a Brauer factor set there is a \(2\)-cocycle in our notation.

\begin{definition}\label{def:BrauerCSA}
    Let \(c \in Z'^2(k,F)\). Let \(B(F,c)\) be the \(k\)-vector space of \(G\)-homogeneous maps from \(\Phi^2\) to \(K\). We turn \(B(F,c)\) into a \(k\)-algebra using the following definition. If \(\ell,\ell' \in B(F,c)\), we set \(\ell \ell' = \ell''\) such that for all \(i,j \in [d]\),
    \[\ell''_{i,j} = \sum_{k=1}^d \ell_{i,k} c_{i,k,j} \ell'_{k,j}.\]
\end{definition}

\begin{fact}\label{thm:BrauerFSAssociated}
    If \(c\) is a \(2\)-cocycle, the \(k\)-algebra \(B(F,c)\) is a central simple \(k\)-algebra. Furthermore, if \(c_1\) and \(c_2\) are associated cocycles, and \(a\) is a \(1\)-cochain such that \(c_2 = \Delta^1(a) c_1\), then the map
    \[\begin{array}{ccl} B(F,c_1) &\to &B(F,c_2) \\ (m_{i,j})_{i,j \in [d]} &\mapsto &(m_{i,j}a_{i,j}^{-1})_{i,j \in [d]}\end{array}\]
    is an isomorphism. Conversely, if \(A\) is a central simple \(k\)-algebra split by \(F\), then there exists a cocycle \(c \in Z^2(k,F)\) such that \(A \simeq B(F,c)\) If \(c_1,c_2 \in Z'^2(k,F)\) are such that \(B(F,c_1)\) and \(B(F,c_2)\) are isomorphic algebras, then the cocycles \(c_1\) and \(c_2\) are associated.
\end{fact}

\begin{proof}
    This is essentially the content of \cite[Theorem 2.5.6]{jacobson1996finite}. While the result is stated in the source for reduced factor sets, we note that the isomorphism constructed in \cite[Theorem 2.3.21]{jacobson1996finite} is valid even for non-reduced associated cocycles. Then, the results stated here extend to non-reduced cocycle by \Cref{prop:ReducedUseless}.
\end{proof}

\begin{remark}\label{rem:ReducedAssociated}
    We note that if \(c_1, c_2 \in Z'^n(k,F)\) are associated and both reduced, then a cochain \(a \in C^{n-1}(k,F)\) such that \(c_1 = \Delta^{n-1}(a) c_2\) must be reduced as well.
\end{remark}

We refer the reader to \cite[Section 2.3]{jacobson1996finite} for the detailed construction of a Brauer factor-set associated to a given central simple algebra, and we give a brief description below:
Let \(A\) be a central simple \(k\)-algebra, and let \(F\) be a separable maximal commutative subalgebra of \(A\). Then, there exists \(v \in A\) such that \(A = FvF\). Furthermore, letting \(K\) be a Galois extension of \(k\) which splits \(F\), there is a map \(A \to M_d(K)\) such that every \(u \in F\) is mapped to \(\Diag(\phi_1(x),\hdots,\phi_d(x))\). Let \((v_{ij})_{i,j \in [d]}\) be the image of \(v\) in \(M_d(K)\). Then, every element \(uvu'\) of \(A\) becomes a matrix of the form \((\ell_{ij} v_{ij})_{i,j \in [d]}\), where \(\ell_{ij} = \phi_i(u)\phi_j(u')\). The map \(uvu' \mapsto (\ell_{ij})_{i,j \in [d]}\) gives an isomorphism from \(A\) to \(B(F,c)\), where \(c_{ijk} \coloneqq v_{ik}v_{kj}v_{ij}^{-1}\).


\section{Explicit computations with Amitsur cohomology}\label{Sec:Amitsur}
In this section, we recall the basic definitions of Amitsur cohomology and then use them to get efficient representations of cocycles as introduced in \Cref{Sec:BrauerPrelim}. We also keep notations as they were in \Cref{Sec:BrauerPrelim}.

Although the Amitsur complex and its cohomology are defined for general commutative rings, we will only need to state the definitions over fields and étale algebras. We refer to \cite[Chapter 5]{ford2017separable} for more general definitions and results.

Unless specified otherwise, all tensor products are taken over \(k\), and by \(F^{\otimes n}\) we mean the tensor product of \(n\) copies of the \(k\)-algebra \(F\).

\begin{definition}
    An \(n\)-cochain in the sense of Amitsur is an invertible element in the \(k\)-algebra \(F^{\otimes n+1}\), and we write \(C_{Am}^n(k,F)\) for the group of \(n\)-cochains. 

    For \(n \in \N\) and \(0 \leq i \leq n+1\), we define a \(k\)-algebra homomorphism \(\varepsilon_i^n\) from \(F^{\otimes n+1}\) to \(F^{\otimes n+2}\). The map \(\varepsilon_i^n\) is defined on the simple tensor elements as follows:
    \[\varepsilon_i^n(f_0 \otimes \hdots \otimes f_n) = f_0 \otimes \hdots \otimes f_{i-1} \otimes 1 \otimes f_{i} \otimes \hdots \otimes f_n.\]
    We then define the group homomorphism
    \[\Delta_{Am}^n\colon\begin{array}{ccl} C_{Am}^n(k,F) &\to &C_{Am}^{n+1}(k,F) \\ a &\mapsto &\prod_{i=0}^{n+1} \varepsilon_i^n(a)^{(-1)^i}\end{array}.\]

    Define the subgroup \(Z_{Am}^n(k,F) = \Ker(\Delta_{Am}^n)\) of \(C_{Am}^n(k,F)\), and its elements are called \(n\)-cocycles in the sense of Amitsur.

    If \(n\) is positive, we also define the subgroup \(B_{Am}^n(k,F) = \im(\Delta_{Am}^{n-1})\) of \(Z_{Am}^n(k,F)\) and its elements are \(n\)-coboundaries in the sense of Amitsur. Then, the \(n\)th Amitsur cohomology group \(H_{Am}^n(k,F)\) is defined as the quotient group \(Z_{Am}^n(k,F)/B_{Am}^n(k,F)\).

    Two \(n\)-cocycles \(c_1\) and \(c_2\) are called associated if they have the same class in \(H_{Am}^n(k,F)\). If \(b \in B^n_{Am}(k,F)\), a cochain \(a \in C^{n-1}_{Am}(k,F)\) such that \(\Delta_{Am}^{n-1}(a) = c\) is called a \emph{trivialisation} of \(c\).
\end{definition}

\subsection{Amitsur cohomology and central simple algebras}\label{Sec:AmitsurAlgebra}

    The first result we need to leverage Amitsur cohomology in representing Brauer factor sets is the following
    \begin{lemma}\label{lemma:AdamsonAmitsur}
        Consider the map \(\Psi_n\) from \(F^{\otimes n+1}\) to the algebra of \(G\)-homogeneous maps from \(\Phi^{n+1}\) to \(K\), which is defined over simple tensors by
        \[f_0 \otimes \hdots \otimes f_n \mapsto \left(\begin{array}{ccl} \Phi^{n+1} &\to &K\\ (\varphi_{i_0},\hdots,\varphi_{i_n}) &\mapsto &\prod_{j=0}^n \varphi_{i_j}(f_j)\end{array}\right).\]
        The map \(\Psi\) is an isomorphism of \(k\)-algebras, which also sends the unit group \(C_{Am}^n(k,F)\) onto the unit group \(C^n(k,F)\). Furthermore, \(\Psi_{n+1} \circ \Delta_{Am}^n = \Delta^n \circ \Psi_{n}\), and it follows that cocycle and coboundaries subgroup are also isomorphic to one another, and so are cohomology groups.
    \end{lemma}

    \begin{proof}
        If \(F\) is a separable field extension of \(k\), this is the content of \cite[Section 2]{rosenberg1960amitsur}. However, the proof given there readily generalises to the case that \(F\) is an étale \(k\)-algebra.
    \end{proof}

    If \((\varphi_{i_0},\hdots,\varphi_{i_n}) \in \Phi^{n+1}\), we get a homomorphism of \(k\)-algebras
        \[\begin{array}{ccc}
            F^{\otimes n+1} &\to &K \\
            a_0 \otimes \hdots \otimes a_n &\mapsto &\prod_{j=0}^n \varphi_{i_j}(a_j)
        \end{array}\]
    By the construction of \(\Psi\), if \(a \in F^{\otimes n+1}\) and \((\varphi_{i_0},\hdots,\varphi_{i_n}) \in \Phi^{n+1}\), then we get
    \begin{equation}\label{eq:AmitsurEmbeddings}
        \Psi_n(a)_{i_0,\hdots,i_n} = (\varphi_{i_0},\hdots,\varphi_{i_n})(a).
    \end{equation}

    Adapting \Cref{def:BrauerCSA} to Amitsur cohomology, we get:
    \begin{definition}\label{def:AmitsurCSA}
        Let \(c \in Z_{Am}^2(k,F)\) be a cocycle. The \(k\)-algebra \(A(F,c)\) is defined as follows:

        As a \(k\)-vector space, \(A(F,c)\) is \(F^{\otimes 2}\). We see \(F^{\otimes 3}\) as a \(F^{\otimes 2}\)-algebra via the embedding \(\varepsilon_1^1\colon F^{\otimes 2} \to F^{\otimes 3}\). Then, multiplication in \(A(F,c)\) is defined by
        \[a a' = \Tr_{F^{\otimes 3}/F^{\otimes 2}}\left(\varepsilon_2^1(a)c\varepsilon_0^1(a')\right).\]
        In the sequel, we use the same embedding whenever we write \(\Tr_{F^{\otimes 3}/F^{\otimes 2}}\).
    \end{definition}

    Observe that the algebra \(F^{\otimes 3}/F^{\otimes 2}\) may be seen as the extension of scalars to \(F^{\otimes 2}\) of the \(k\)-algebra \(F\). Since we make the choice of seeing \(F^{\otimes 3}\) as a \(F^{\otimes 2}\)-algebra via the mapping \(\varepsilon_1^2\), the relevant injection \(F \to F^{\otimes 3}\) sends an element \(a \in F\) to \(1 \otimes a \otimes 1\). Then, by \cite[Section III.9.3]{bourbaki1998algebra}, for \(a \in F\), we have
    \[\Tr_{F^{\otimes 3}/F^{\otimes 2}}(1 \otimes a \otimes 1) = \Tr_{F/k}(a) (1 \otimes 1).\]
    By linearity of the trace, if \(a = \sum_{i \in I} a_{i0} \otimes a_{i1} \otimes a_{i2} \in F^{\otimes 3}\), we get
        \[\Tr_{F^{\otimes 3}/F^{\otimes 2}}(a) = \sum_{i \in I} \Tr_{F/k}(a_{i1}) (a_{i0} \otimes a_{i2}).\]

    \begin{prop}\label{prop:IsomBrauerAmitsur}
        Let \(c \in Z_{Am}^2(k,F)\). The map \(\Psi_1\) is a \(k\)-algebra isomorphism from \(A(F,c)\) to \(B(F,\Psi_2(c))\).
    \end{prop}

    \begin{proof}
        The map \(\Psi_1\) is already known to be an isomorphism of vector space between \(A(F,c)\) and \(B(F,\Psi_2(c))\), as the underlying vector spaces of these algebras are respectively \(F^{\otimes 2}\) and the space of \(G\)-homogeneous maps from \(\Phi^2\) to \(K\). We therefore only need to check that the map \(\Psi_1\) is a homomorphism of \(k\)-algebras.

        It is well known that if \(a \in F\), \(\Tr_{F/k}(a) = \sum_{\ell \in [d]}\varphi_\ell(a)\). Then, if \(a = a_0 \otimes a_1 \otimes a_2 \in F^{\otimes 3}\) and \(i,j \in [d]\), we get:
        \begin{align*}
        (\varphi_i,\varphi_j)\left(\Tr_{F^{\otimes 3}/F^{\otimes 2}}(a)\right) &= (\varphi_i,\varphi_j)\left(\Tr_{F/k}(a_1) (a_0 \otimes a_2)\right) \\
        &= (\varphi_i,\varphi_j)\left(\sum_{\ell \in [d]} \varphi_\ell(a_1) (a_0 \otimes a_2)\right) \\
        &= \sum_{\ell \in [d]} \varphi_i(a_0) \varphi_\ell(a_1) \varphi_j(a_2) \\
        &= \sum_{\ell \in [d]} (\varphi_i,\varphi_\ell,\varphi_j)(a_0 \otimes a_1 \otimes a_2).
        \end{align*}
        Extending this result by linearity, we get that 
        \[(\varphi_i,\varphi_j) \circ \Tr_{F^{\otimes 3}/F^{\otimes 2}} = \sum_{\ell \in [d]} (\varphi_i,\varphi_\ell,\varphi_j).\]

        We also note that if \(a \in F^{\otimes 2}\), then by \Cref{eq:AmitsurEmbeddings} we get
            \[(\varphi_i,\varphi_\ell,\varphi_j)(\varepsilon_2^1(a)) = (\varphi_i,\varphi_\ell)(a) = \Psi_1(a)_{i,\ell}.\]
            Likewise,
            \[(\varphi_i,\varphi_\ell,\varphi_j)(\varepsilon_0^1(a)) = \Psi_1(a)_{\ell,j}.\]

        Now, using \Cref{eq:AmitsurEmbeddings} again, we get that if \(a,a' \in A(F,c)\) and \(i,j \in [d]\),
        \begin{align*}
            \Psi_1(aa')_{i,j} &= (\varphi_i,\varphi_j)\left( \Tr_{F^{\otimes 3}/F^{\otimes 2}}\left(\varepsilon_2^1(a)c\varepsilon_0^1(a')\right)\right)\\
            &= \sum_{\ell \in [d]}(\varphi_i,\varphi_\ell,\varphi_j)\left(\varepsilon_2^1(a)c\varepsilon_0^1(a')\right) \\
            &= \sum_{\ell \in [d]}(\varphi_i,\varphi_\ell,\varphi_j)(\varepsilon_2^1(a))(\varphi_i,\varphi_\ell,\varphi_j)(c)(\varphi_i,\varphi_\ell,\varphi_j)(\varepsilon_0^1(a')) \\
            &= \sum_{\ell \in [d]} \Psi_1(a)_{i,\ell} \Psi_2(c)_{i,\ell,j} \Psi_1(a')_{\ell,j} \\
            &= \Psi_1(a) \Psi_1(a'),
        \end{align*}
        where the multiplication in the last line is meant to be in \(B(F,\Psi_2(c))\).
    \end{proof}

    \begin{remark}\label{remark:AmitsurTrivialCocycleIsom}
        Let \(\1 \in F^{\otimes 3}\) be the trivial cocycle. There is an explicit isomorphism between \(A(F,\1)\) and \(M_d(k)\). Indeed, the algebra \(M_d(k)\) is isomorphic to the algebra \(F \otimes F^\vee\) of \(k\)-linear endomorphisms of \(F\), where \(F^\vee\) is the dual \(k\)-vector space of \(F\). In fact, the map \(a \mapsto \Tr_{F/k}(a \cdot) \coloneq \left(b \mapsto \Tr_{F/k}(ab)\right)\) gives a \(k\)-linear isomorphism from \(F\) to its dual space \(F^\vee\), so we may see \(M_d(k)\) as the algebra \(F \otimes F\) with the following multiplication law:
        \[(a \otimes a')(b \otimes b') = \Tr_{F/k}(a' b) a \otimes b'.\]
        Let \(a,a',b,b' \in F\), we compute the following multiplication in \(A(F,\1)\):
        \begin{align*}
            (a \otimes a') (b \otimes b') &= \Tr_{F^{\otimes 3}/F^{\otimes 2}}\left(\varepsilon_2^1(a \otimes a') \varepsilon_0^1(b \otimes b')\right) \\
                                            &= \Tr_{F^{\otimes 3}/F^{\otimes 2}}(a \otimes a'b \otimes b') \\
                                            &= \Tr_{F/k} (a' b) a \otimes b'.
        \end{align*}

        It follows that if \(\mathcal{B}\) is a basis of \(F\) and we identify elements of \(F\) with the column vectors of their components in basis \(\mathcal{B}\), and for \(a \in F\), we set \(\varphi(a)\) to be the row matrix of the linear form \(b \mapsto \Tr_{F/k}(ab)\), then an isomorphism from \(A(F,\1)\) to \(M_d(k)\) is given by
        \[a \otimes a' \mapsto a\varphi(a').\]
    \end{remark}
    
    \begin{theorem}\label{thm:IsomAmiAlgebra}
        Let \(A\) be a central simple \(k\)-algebra, let \(F \subset A\) be a separable maximal commutative subalgebra, and let \(v \in A\) be such that \(A = FvF\). Then, there is a cocycle \(c \in Z^2_{Am}(k,F)\) such that the map
        \[\begin{array}{rlcl}
            \Phi_v\colon &F^{\otimes 2} &\to &A \\
            &r_1 \otimes r_2 &\mapsto &r_1 v r_2
        \end{array}\]
        is an isomorphism of \(k\)-algebras from \(A(F,c)\) to \(A\) (which makes sense as \(F^{\otimes 2}\) is the underlying \(k\)-vector space of \(A(F,c)\)). Furthermore, \(c\) is the unique element of \(F^{\otimes 3}\) with this property (if we extend the definition of the multiplication in \(A(F,c)\) to non-cocycle elements).
    \end{theorem}

    \begin{proof}
        Let \(c' \in B^2(k,F)\) be the Brauer factor set constructed from the data of \(F, A\) and \(v\) as in the discussion at the end of \Cref{Sec:BrauerPrelim}, and let \(\theta_v\) be the isomorphism from \(A\) to \(B(F,c')\) constructed there. We also let \(c = \Psi_2^{-1}(c')\). The fact that \(\Phi_v\) is an isomorphism of \(k\)-algebras follows from the observation that \(\theta_v \circ \Phi_v = \Psi_1\) and from \Cref{prop:IsomBrauerAmitsur}.
        
        The cocycle \(c\) is a solution to the system of \(F^{\otimes 2}\)-linear equations
        \[\Tr_{F^{\otimes 3}/F^{\otimes 2}} \left(\varepsilon_2^1(u)c\varepsilon_0^1(u')\right) = \Phi_v^{-1}\left(\Phi_v(u)\Phi_v(u')\right)\]
        for all \(u,u' \in F^{\otimes 2}\). Observe that \(\varepsilon_2^1(u) \varepsilon_0^1(u')\) ranges over \(F^{\otimes 3}\) as \(u,u'\) range over \(F^{\otimes 2}\). Then, the uniqueness of \(c\) follows from the fact that the map \(c \mapsto \Tr_{F^{\otimes 3}/F^{\otimes 2}}(c \cdot )\) is an isomorphism from \(F^{\otimes 3}\) to its dual as a \(F^{\otimes 2}\)-module.

        This last fact, in turn, follows from \cite[Corollary 4.6.8 and Proposition 4.6.1]{ford2017separable}. Indeed, \(F^{\otimes 3} \simeq F^{\otimes 2}[X]/P\), and \(P \in k[X] \subset F^{\otimes 2}[X]\) is a separable polynomial, so \(F^{\otimes 3}\) is a separable \(F^{\otimes 2}\)-algebra.
    \end{proof}

    \begin{cor}\label{cor:AssoCocyIsomAmitsur}
        Let \(c,c' \in Z_{Am}^2(k,F)\) be associated cocycles, and let \(a \in C_{am}^2(k,F)\) be such that \(c' = c\Delta_{Am}^1(a)\). Then, the map \(m \mapsto ma^{-1}\), where the multiplication used is the natural one in \(F^{\otimes 2}\), is an isomorphism from \(A(F,c)\) to \(A(F,c')\).
    \end{cor}

    \begin{proof}
        This is just a consequence of the analogous property for Brauer factor sets (\Cref{thm:BrauerFSAssociated}), pulled back through the isomorphisms \(\Psi_i\).
    \end{proof}

    \begin{remark}
        For the sake of efficiency, we use already proven results on Brauer factor-sets and the isomorphisms \(\Psi_1\) and \(\Psi_2\) to describe the isomorphism between the groups \(H^2_{Am}(k,F)\) and \(\Br(F/k)\). However, the simplicity of the isomorphism \(\Phi_c\) described in \Cref{thm:IsomAmiAlgebra} suggests that Amitsur cohomology is a natural setting for giving an algebraic description of \(\Br(F/k)\) and that direct proofs of these facts may be given without using references to Brauer factor-sets.
    \end{remark}

\subsection{Computational representation of Amitsur cohomology}\label{Sec:CompRepAmi}
    Once a field \(k\) and an étale \(k\)-algebra \(F = k[X]/P\) are fixed, for \(n \in \N\), we have an isomorphism \[\Xi_n\colon F^{\otimes n+1} \to R_n \coloneqq k[X_0,\hdots,X_n]/I_n,\] where \(I_n\) is the ideal generated by the polynomials \(P(X_i)\) for \(0 \leq i \leq n\). Furthermore, an element of \(R_n\) admits a unique lift in \(k[X_0,\hdots,X_n]\) such that the individual degrees in \(X_0,X_1,\hdots,X_n\) are all lesser or equal to \(d-1\). If \(Q \in k[X_0,\hdots,X_n]\), we denote by \(\widetilde{Q}\) the lift of \(Q \mod I\) as described above. There is an obvious polynomial algorithm which, given some \(Q \in k[X_0,\hdots,X_n]\), computes \(\widetilde{Q}\). Indeed, one may simply reduce \(Q\) modulo \(P(X_0)\), then modulo \(P(X_1)\), and so on... This allows us to represent an element \(Q \in R_n\) as its lift \(\widetilde{Q}\), and in general we get \(Q_1 + Q_2 = \widetilde{Q_1} + \widetilde{Q_2} \mod I_n\) and \(Q_1 Q_2 = \widetilde{Q_1} \widetilde{Q_2} \mod I_n\).

    In this setting, the map \(\varepsilon_i^n\) becomes the map \[Q \mod I_n \mapsto \widetilde{Q}(X_0,\hdots,X_{i-1},X_{i+1},\hdots,X_n) \mod I_{n+1}.\] That is, we define a map
    \[\nu_i^n(X_j) = \begin{cases}
        X_j \quad \text{if } j < i \\
        X_{j+1} \quad \text{otherwise.}
    \end{cases}
    \]
    from \(k[X_0,\hdots,X_n]\) to \(k[X_0,\hdots,X_{n+1}]\), and we define
    \[\begin{array}{rccc}
        \varepsilon_i^n\colon &R_n &\to &R_{n+1} \\
        &   Q  \mod I_n &\mapsto &\nu_i^n(\widetilde{Q}) \mod I_{n+1}.
    \end{array}\]
    The maps defined above are compatible with the eponymous maps from the Amitsur complex. That is, we have \(\varepsilon_i^n \circ \Xi_n = \Xi_{n+1} \circ \varepsilon_i^n\).

    Now, the trace map of \(R_2\) as an \(R_1\) algebra (via the morphism \(\varepsilon_1^1\)) may easily be computed in the \(R_1\)-basis \((X_1^i)_{0 \leq i < d}\) of \(R_2\). Let \(c \in Z^2_{Am}(k,F)\), \(Q = \Xi_2(c)\), \(a,a' \in A(F,c)\), \(T = \Xi_1(a)\) and \(T' = \Xi_1(a')\). We set \(T'' = \Xi_1^{-1}(aa')\), and we have
    \[T'' = \Tr_{R_2/R_1}\left(\varepsilon_2^1(T)Q\varepsilon_0^1(T')\right).\]
    That is,
    \[\Xi_1(aa') = \Tr_{R_2/R_1}\left(\varepsilon_2^1(\Xi_1(a))\Xi_2(c)\varepsilon_0^1(\Xi_1(a'))\right).\]

    For the remainder of this work, elements of algebras of the form \(F^{\otimes n+1}\), when \(F \simeq k[X]/(P(X))\) is an étale \(k\)-algebra, are represented algorithmically as their images by \(\Xi_n\) in \(R_n\).
\subsection{Representing central simple algebras as Amitsur algebras}\label{Sec:AmitsurAlgebraConstruction}

    We may now give an algorithm for finding a representation of a central simple algebra \(A\) over any field of sufficiently large size where linear algebra tasks may be performed efficiently.

    \begin{algorithm}
        \KwIn{A field \(k\)}
        \KwIn{A central simple \(k\)-algebra \(A\) such that \(|k| > \dim_k A\)}
        \KwOut{\(u \in A\) such that \(F = k(u)\) is a maximal commutative subalgebra of \(A\)}
        \KwOut{\(P \in k[X]\), the minimal polynomial of \(u\)}
        \KwOut{\(c \in Z_{Am}^2(k,F)\)}
        \KwOut{A linear map \(e\colon F^{\otimes 2} \to A\) that is an isomorphism from \(A(F,c)\) to \(A\)}
        Find \(u \in A\) such that \(F \coloneqq k[u]\) is a maximal separable commutative subalgebra of \(A\) \label{algoline:Findu}\;
        Compute \(P\), the minimal polynomial of \(u\)\;
        Find \(v \in A\) such that \(A = FvF\)\; \label{algoline:Findv}
        Compute the \(k\)-vector space isomorphism \(e: F^{\otimes 2} \to A\) sending \(u^i \otimes u^j\) to \(u^i v u^j\)\;
        Compute \(c \in F^{\otimes 3}\) such that for \(0 \leq i,j,i',j' \leq d-1\), \(e(u^i \otimes u^j)e(u^{i'} \otimes u^{j'}) = \Tr_{F^{\otimes 3}/F^{\otimes 2}}(\varepsilon_2^1(u^i \otimes u^j)c\varepsilon_0^1(u^{i'} \otimes u^{j'}))\)\; \label{algoline:Findc}
        \KwRet{\((u,P,c,e)\)}
        \caption{Computing a presentation of a central simple algebra as an Amitsur algebra}
        \label{algo:FindCocycle}
    \end{algorithm}

    Before we prove the correctness and efficiency of \Cref{algo:FindCocycle}, we need a lemma:

    \begin{lemma}\label{lemma:Findv}
        Let \(k\) be a field and let \(A\) be a central simple \(k\)-algebra. Assume that \(|k| > \dim_k A\). Let \(u \in A\) be such that \(F \coloneqq k[u]\) is a maximal commutative subalgebra of \(A\). Then an element \(v \in A\) such that \(A = FvF\) may be found in probabilistic polynomial time.
    \end{lemma}

    \begin{proof}
        For \(v\) in \(A\), by an argument of dimensions over \(k\), we observe that \(A = FvF\) if and only if the map
        \[e\colon\begin{array}{ccl} F \otimes F &\to &A \\ f_1 \otimes f_2 &\mapsto &f_1vf_2 \end{array}\]
        is injective.

        The family \((u^i \otimes u^j)_{0 \leq i,j \leq \deg A - 1}\) is a basis of \(F^{\otimes 2}\). Let \(B = (b_1,\hdots,b_{\dim A})\) be the input basis of \(A\) (that is, the basis for which the structure constants of \(A\) are expressed). Then, identifying \(e\) with its matrix in the bases written above, the determinant of \(e\) is a homogeneous polynomial on the coordinates of \(v\) in basis \(B\) of degree bounded by \(\dim A\). Furthermore, \(A = FvF\) if and only if \(v\) is not a zero of this polynomial. We note that by \cite[Theorem 2.2.2]{jacobson1996finite}, this polynomial is nonzero. 

        Letting \(S\) be a finite subset of \(k\), the Schwartz-Zippel lemma ensures that a random \(v\) in \(Sb_1 \oplus \hdots \oplus Sb_{\dim A}\) satisfies this condition with probability larger than \(1 - \frac{\dim A}{|S|}\).

        Therefore, if \(|k| > \dim A\), we may pick \(S\) large enough that \(v\) has the desired property with positive probability, and small enough that we may draw a random element of \(S\) efficiently. For instance, take \(|S| = \dim A + 1\). Then, \(v\) has the desired property with probability larger than \(\frac{1}{\dim A + 1}\). In this case, it is expected to take less than \(\dim A + 1\) attempts to find a suitable element \(v \in A\).
    \end{proof}
    
    \begin{theorem}\label{thm:AlgoFindCocycle}
        If \(k\) is a field over which linear algebra may be performed efficiently, and \(A\) is a central simple \(k\)-algebra such that \(|k| > \dim_k A\), then \Cref{algo:FindCocycle} is correct and runs in probabilistic polynomial time.
    \end{theorem}

    \begin{proof}
         We first prove the correctness of the algorithm. The first 3 steps compute \(u,v \in A\), as well as \(P \in k[X]\), the minimal polynomial of \(u\). We obtain a subalgebra \(F = k[u] \simeq k[X]/P\) of \(A\) and \(v \in A\) such that \(A = FvF\). Then, by \Cref{thm:IsomAmiAlgebra}, there exists a cocycle \(c \in Z_{Am}^2(k,F)\) such that \(e\colon a \otimes b \mapsto avb\) is an isomorphism from \(A(k,F)\) to \(A\), and such a \(c\) is unique in \(F^{\otimes 3}\). The equation solved in \Cref{algoline:Findc} finds a representation of this \(c\).

        Now we consider the complexity of the algorithm:
        
        The element \(u \in A\) in \Cref{algoline:Findu} may be found using the polynomial algorithm given in \cite{graaf2000finding}.

        Then, \(v\) in \Cref{algoline:Findv} may be found in probabilistic polynomial time using the algorithm of \Cref{lemma:Findv}.

        Then, the remaining lines involve arithmetic in the \(k\)-algebras \(A\), \(R_1\) and \(R_2\), as well as the computation of the solution of a system of linear equations.
        
        All in all, this makes \Cref{algo:FindCocycle} a polynomial probabilistic algorithm.
    \end{proof}

\section{Trivialisation of Amitsur cocycles}\label{Sec:Trivialisation}
In this section, we present an algorithm for computing a trivialisation of a coboundary using \(S\)-unit group computation. For this purpose, we prove \Cref{thm:SUnitTriv}, which states that a trivialisation of a coboundary may be found in an appropriate group of \(S\)-units. This result is to be compared with \Cref{thm:Simon} and \cite[Theorem 7]{fieker2009minimizing}. Then, we may state \Cref{algo:TrivCobound} which computes a trivialisation of a given coboundary by solving the related equation in an appropriate group of \(S\)-units. Furthermore, we get \Cref{cor:QuantumSplit} which, assuming GRH, solves the explicit isomorphism problem in quantum polynomial time.

The technical part of this section is devoted to proving \Cref{thm:SUnitTriv}. The strategy we employ is similar to the proof of \cite[Theorem 7]{fieker2009minimizing}. That is, we prove \Cref{lemma:Hilbert90}, which is a generalisation to divisors of Hilbert's Theorem 90 on the vanishing of the first cohomology group. While \cite[Lemma 9]{fieker2009minimizing} is stated for fractional ideals of a number field, our setting requires us bring this machinery to the setting of étale algebras. We therefore first introduce definitions and basic results for the divisors of an étale algebra over a global field in \Cref{sec:DivisorsEtaleAlgebras}.

\subsection{Divisors of étale algebras}\label{sec:DivisorsEtaleAlgebras}
    
    Let \(k\) be a global field. Then, recall that a commutative \(k\)-algebra \(F\) is \emph{étale} if and only if it satisfies the following equivalent conditions:
    \begin{enumerate}
        \item The algebra \(F\) factors as a direct product \(F \simeq F_1 \times \hdots \times F_r\) of finite separable extensions of \(k\).
        \item The algebra \(F\) is isomorphic to an algebra of the form \(k[X]/(P(X))\), where \(P \in k[X]\) is a separable polynomial (that is, \(P\) has no multiple root in an algebraic closure of \(k\)).
    \end{enumerate}

    \begin{remark}
        If \(F\) is an étale \(k\)-algebra, the isomorphism \(F \simeq F_1 \times \hdots \times F_r\) discussed above is unique up to reindexing and automorphism of the factors. In fact, the factors identify with the minimal ideals of \(F\). Therefore, we usually identify an étale algebra \(F\) with the product \(F_1 \times \hdots \times F_r\). We will also denote by \(\pi_i\) the projection map \(F \to F_i\).
    \end{remark}

    If \(F\) is a global field, we write \(\Pl(F)\) for the set of non-archimedean places of \(F\). Then, recall that the divisor group \(\Dc(F)\) is the free abelian group on the set \(\Pl(F)\), that an element \(a \in F\) has an associated divisor \(\Dc(a) = \sum_{P \in \Pl(F)} \ord_P(a) P\), and that if \(\Pc(F) = \{\Dc(a)\colon a \in F\}\) is the group of principal divisors, then the class group \(\cl(F) = \Dc(F)/\Pc(F)\) of \(F\) is finite. When \(F\) is a number field, the divisor group is isomorphic to the group of fractional ideals of \(F\), and the class group in terms of divisors is then isomorphic to the ideal class group.

    \begin{definition}
        Let \(F = F_1 \times \hdots \times F_r\) be an étale \(k\)-algebra factored as above. The \emph{set of non-archimedean places} of \(F\) is the disjoint union \(\Pl(F) = \bigsqcup_{i \in [r]} \Pl(F_i)\). The \emph{divisor group} \(\Dc(F)\) of \(F\) is the free abelian group on \(\Pl(F)\). The divisor \(\Dc(a)\) of an invertible element \(a \in F^\times\) is defined as the sum of the divisors of its projections in the factors \(F_i\) of \(F\), and the groups \(\Pc(F)\) and \(\cl(F)\) are defined by analogy with their definitions for global fields. If \(D \in \Dc(F) = \sum_{P \in \Pl(F)} n_P P\), we let \(\Supp(D) = \{P \in \Pl(F)\colon n_P \neq 0\}\) be the \emph{support} of \(D\).
    \end{definition}

    An immediate consequence of this definition is that the class group of an étale \(k\)-algebra is finite.

    If \(F = F_1 \times \hdots \times F_r\) and \(F' = F'_1 \times \hdots \times F'_{r'}\) are two étale \(k\)-algebras, and \(f\colon F \to F'\) is a homomorphism of \(k\)-algebras, it decomposes as a matrix \((f_{ij})_{\substack{i \in [r] \\ j \in [r']}}\), where \(f_{ij}\colon F_i \to F_j'\), and each \(f_{ij}\) is either the zero map or an embedding of field extensions of \(k\) (depending on whether the unit of \(F_i\) is mapped to \(0\) or to the unit of \(F_j\)). That is, if \(a \in F\) and \(j \in [r']\),
    \[\pi'_j(f(a)) = \sum_{i \in [r]} f_{ij}(\pi_i(a)),\]
    where the \(\pi_i\) are the projections \(F \to F_i\) and the \(\pi'_j\) are the projections \(F' \to F'_j\).

    We may then define a map \(\Dc(f)\colon \Dc(F) \to \Dc(F')\) the following way: let \(i \in [r]\) and \(P \in \Pl(F_i)\), we set \(\Dc(f)(P) = \sum_{j \in [r']} \Dc(f_{ij})(P)\), where \(\Dc(f_{ij})\) is the usual map on divisors of global fields if \(f_{ij}\) is an embedding of fields and is zero if \(f_{ij} = 0\). One may check that this definition makes \(\Dc\) into a functor from the category of étale \(k\)-algebras to the category of abelian groups. Furthermore, the map \(\Dc\colon F^\times \to \Dc(F)\) is a natural transformation of the multiplicative group functor into the functor \(\Dc\) sending an étale algebra to its divisor group. That is, if \(a \in F^\times\), then \(\Dc(f)(\Dc(a)) = \Dc(f(a))\).

    As is the case for global fields, if \(Q \in M_{F'}\), there is at most one place \(P \in \Pl(F)\) such that \(Q \in \Supp(\Dc(f)(P))\). To prove this, we first need a lemma:
    \begin{lemma}\label{lemma:UniqueFieldPreimage}
        Let notations be as above and let \(j \in [r']\). Then there is at most one \(i \in [r]\) such that \(f_{ij}\) is nonzero.
    \end{lemma}

    \begin{proof}
        Let \(i_1,i_2 \in [r]\) and \(j \in [r']\), and let \(e_1,e_2\) be the units respectively of \(F_{i_1}\) and \(F_{i_2}\). We then have \(\pi'_j(f(e_1)) = f_{i_1j}(1) \in F_j\) and \(\pi'_j(f(e_2)) = f_{i_2j}(1) \in F_j\). Now, since \(e_1 e_2 = 0\), we have \(f_{i_1j}(1) f_{i_2j}(1) = 0\). That is, either \(f_{i_1j}(1) = 0\) or \(f_{i_2j}(1) = 0\). It follows that either \(f_{i_1j}\) or \(f_{i_2j}\) is the zero map.
    \end{proof}

    Now, we prove the following:
    \begin{prop}\label{prop:UniquePlaceAbove}
        Let notations be as above, and let \(Q \in \Pl(F')\). Then there is at most one place \(P \in \Pl(F)\) such that \(Q \in \Supp(\Dc(f)(P))\).
    \end{prop}

    \begin{proof}
        Let \(j \in [r']\) be such that \(Q \in \Pl(F_j)\). Then, by \Cref{lemma:UniqueFieldPreimage}, there is at most one value \(i \in [r]\) such that \(f_{ij} \neq 0\). If there is no such \(i\), then no divisor of \(F\) has an image by \(\Dc(f)\) with \(Q\) in its support. Otherwise, if \(Q \in \Supp(\Dc(f)(P))\), it follows that \(P \in \Pl(F_i)\). Then, the usual theory of divisors of global fields ensures that there is only one \(P \in \Pl(F_i)\) such that \(Q \in \Supp(\Dc(f)(P))\).
    \end{proof}

    If \(j \in [r']\) is such that \(f_{ij}\) is nonzero for some \(i \in [r]\), and if \(Q \in \Pl(F'_j)\), we write \(Q_f\) for the unique \(P \in \Pl(F)\) such that \(Q \in \Supp(\Dc(f)(P))\). We note that if \(f\) is an inclusion of the form \(F \hookrightarrow F \otimes_k F'\), \(a \mapsto a \otimes 1\), then \(Q_f\) exists for all \(Q \in \Pl(F \otimes_k F')\).

    Let \(Q \in \Pl(F')\) such that \(Q_f\) is well defined. We let \(e_{Q,f}\) be the coefficient of \(Q\) in the expansion of the divisor \(\Dc(f)(Q_f)\). Let \(P \in \Pl(F)\), if \(e_{Q,f} \neq 1\) for some \(Q \in \Supp(\Dc(f)(P))\), we say that the place \(P\) \emph{ramifies through \(f\)}. This is equivalent to saying that, if \(i \in [r]\) is such that \(P \in \Pl(F_i)\), then for some \(j \in [r']\) such that \(f_{ij}\) is nonzero, the place \(P\) ramifies in \(F_j\).

    In order to prove some facts about the behaviour of places in étale algebras, we will need to use local fields. For this purpose, we make the following definition:

    \begin{definition}
        Let \(F\) be as above. Let \(i \in [r]\) and let \(P \in \Pl(F_i) \subset \Pl(F)\). Then the \emph{local completion} \(F_P\) of \(F\) at \(P\) is the completion of the global field \(F_i\) at the place \(P\). This completion is an \(F\)-algebra via the composite map
        \[F \xrightarrow{\pi_i} F_i \rightarrow F_P.\]
    \end{definition}

    As in the case of extensions of global fields, the places of \(F'\) lying above a place of \(F\) may be classified as the direct factors of a tensor product.

    \begin{prop}\label{prop:EtaleAlgebraPlaceTP}
        Let notations be as above, let \(i \in [r]\) and \(P \in \Pl(F_i) \subset \Pl(F)\). Let \(J = \{j \in [r']\colon f_{ij} \neq 0\}\). Then, we have
        \[F_P \otimes_F F' \simeq F_P \otimes_{F_i} \prod_{j \in J} F'_j.\]
        The usual bijections between the places of \(F'_j\) above \(P\) and the direct factors of \(F_P \otimes F'_j\) give a bijection between the direct factors of \(F_P \otimes F'\) and \(\Supp(\Dc(f)(P))\).
        Furthermore, for \(Q \in \Supp(\Dc(f)(P)\), we have \(e_{Q,f} \neq 1\) if and only if the corresponding factor is a ramified extension of \(F_P\).
    \end{prop}

    \begin{proof}
        We have
        \[F_P \otimes_{F_i} \prod_{j \in J} F'_j \simeq F_P \otimes_{F_i} (F_i \otimes_F F') \simeq F_P \otimes_F F'.\]
        Then, the rest of the results follows directly from the definition of \(\Dc(f)\) and the equivalent results for extensions of global fields (see \textit{e.g} \cite[Sections II.10 and II.19]{cassels1967algebraic}).
    \end{proof}

    Finally, we also prove a result about ramification and tensor products:
    \begin{prop}\label{prop:TPUnramified}
        Let \(F\) and \(G\) be étale \(k\)-algebras, and let \(P \in \Pl(k)\) be a place that ramifies neither in \(F\) nor in \(G\). Then, \(P\) does not ramify in \(F \otimes_k G\) either.
    \end{prop}

    \begin{proof}
        Let \(k_P\) be the completion of \(k\) at \(P\). By \Cref{prop:EtaleAlgebraPlaceTP}, it is enough to prove that the direct factors of \(k_P \otimes_k F \otimes_k F'\) are unramified extensions of \(k_P\). Now, since
        \[k_P \otimes_k F \otimes_k F' \simeq (k_P \otimes_k F) \otimes_{k_P} (k_P \otimes_k F'),\]
        the result will follow from the fact that if \(K,K'\) are unramified finite extensions of \(k_P\), then the direct factors of \(K \otimes_{k_P} K'\) are themselves unramified extensions of \(k_P\).

        This, in turns, follows from \cite[Proposition I.7.1]{cassels1967algebraic}. Indeed, by the proposition, we know that there are polynomials \(g,g' \in R[X]\) (where \(R\) is the valuation ring of \(k_P\)) such that \(K \simeq k_P[X]/(g(X))\), \(K' \simeq k_P[X]/g'(X)\) and whose residues \(\overline{g}\) and \(\overline{g'}\) are irreducible and separable in \(\kappa_P[X]\), where \(\kappa_P\) is the residue field of \(k_P\). 
        Now, we have \(K \otimes_{k_P} K' \simeq K[X]/(g'(X))\). A direct factor of this algebra is of the form \(K[X]/(h(X))\), where \(h\) is an irreducible factor of \(g'\) in \(\Oc[X]\), where \(\Oc\) is the valuation ring of \(K\). Since \(\overline{h}\) is a factor of the separable polynomial \(\overline{g'}\), it is itself separable as a polynomial in \(\mathfrak{k}[X]\), where \(\mathfrak{k}\) is the residue field of \(K\). Since the polynomial \(h\) is irreducible in \(\Oc[X]\) and its residue is separable, this residue \(\overline{h}\) is irreducible in \(\mathfrak{k}[X]\) by Hensel's lemma. By the proposition cited above, the field extension \(K[X]/(h(X))\) is unramified over \(K\), and therefore over \(k_P\). It follows that the direct factors of \(K \otimes_{k_P} K'\) are unramified over \(k_P\).
    \end{proof}

    \subsection{An algorithm for trivialising Amitsur coboundaries over a number field}

    Let \(k\) be a number field, and let \(F\) be an étale \(k\)-algebra. The Amitsur complex for \(F\) yields a complex
    \[\hdots \to \Dc(F^{\otimes n+1}) \xrightarrow{\Dc(\Delta_{Am}^{n})} \Dc(F^{\otimes n+2}) \to \hdots.\]
of abelian groups, where we define \(\Dc(\Delta_{Am}^n)\) as \(\sum_{0 \leq i \leq n} (-1)^i \Dc(\epsilon^n_i)\).

We give a precise description of the map \(\Dc(\Delta_{Am}^n)\). Let \(Q \in \Pl(F^{\otimes n+2})\). Then, for any \(0 \leq i \leq n+1\), we set \(Q_i = Q_{\varepsilon_i^n}\) and \(e_{Q,i} = e_{Q,\varepsilon_i^n}\). Then, if \(P \in \Pl(F^{\otimes n+1})\), we get
\[\Dc(\epsilon_i^n)(P) = \sum_{\substack{Q \in \Pl(F^{\otimes n+2}) \\ Q_i = P}} e_{Q,i} Q,\]
and it follows that 
\[\Dc(\epsilon_i^n)\left(\sum_{P \in \Pl(F^{\otimes n+1})} n_P P\right) = \sum_{Q \in \Pl(F^{\otimes n+2})} e_{Q,i} n_{Q_i} Q\]
and
\[\Dc(\Delta_{Am}^n)\left(\sum_{P \in \Pl(F^{\otimes n+1})} n_P P\right) = \sum_{Q \in \Pl(F^{\otimes n+2})} \left(\sum_{0 \leq i \leq {n+1}} (-1)^i e_{Q,i} n_{Q_i} \right) Q.\]

We first need two lemmas:

\begin{lemma}\label{lemma:CocycleTransit}
    Let \(Q,Q'\) be places of \(F^{\otimes 2}\) such that \(Q_0 = Q'_0 = P\). Then there exists a place \(R \in \Pl(F^{\otimes 3})\) such that \(R_1 = Q\) and \(R_0 = Q'\).
\end{lemma}

\begin{proof}
    We let \(\alpha\) be a defining polynomial of \(F\). That is, \(F \simeq k[X]/(\alpha(X))\), with \(\alpha \in k[X]\) separable. We may then identify \(F^{\otimes 2}\) with \(F[X]/(\alpha(X))\), where \(\varepsilon_0^0\) maps \(F\) into the ring of scalars in \(F[X]/(\alpha(X))\) and \(\varepsilon_1^0\) is the map
    \[\begin{array}{ccc}
        k[X]/(\alpha(X)) &\to &F[X]/(\alpha(X)) \\
        X &\mapsto &X.
        \end{array}\]

    Likewise, we identify \(F^{\otimes 3}\) with \(F[X,Y]/(\alpha(X),\alpha(Y))\). Then, the descriptions of the maps \(\varepsilon_0^1\) and \(\varepsilon_1^1\) become:
    \[\begin{array}{rccc}
        \varepsilon_0^1\colon &F[X]/(\alpha(X)) &\to &F[X,Y]/(\alpha(X),\alpha(Y)) \\
        & X &\mapsto &Y \\
        \varepsilon_1^1\colon &F[X]/(\alpha(X)) &\to &F[X,Y]/(\alpha(X),\alpha(Y)) \\
        & X &\mapsto &X
        \end{array}\]

    These identifications are coherent with the maps of the Amitsur complex. In particular, we note that \(\varepsilon_0^1 \circ \varepsilon_0^0 = \varepsilon_1^1 \circ \varepsilon_0^0\). With this in place, let \(P \in \Pl(F)\). By \cref{prop:EtaleAlgebraPlaceTP}, the support of \(\Dc(\varepsilon_0^0)(P)\) is in bijection with the direct factors of \(F_P[X]/(\alpha(X))\). That is, it is in bijection with the irreducible factors of \(\alpha(X)\) in \(F_P[X]\). Likewise, the support of \(\Dc(\varepsilon_0^1 \circ \varepsilon_0^0)(P)\) is in bijection with the irreducible factors of \(F_P[X,Y]/(\alpha(X),\alpha(Y))\). That is, it is in bijection with the maximal ideals that contain the ideal \((\alpha(X),\alpha(Y))\).

    Now, if \(Q \in \Supp(\Dc(\varepsilon_0^0)(P))\), we let \(\alpha_Q(X)\) be the irreducible factor of \(\alpha(X)\) corresponding to \(Q\). Likewise, if \(R \in \Supp(\Dc(\varepsilon_0^1 \circ \varepsilon_0^0)(P)\), we let \(\mpp_R\) be the maximal ideal of \(F_P[X,Y]\) corresponding to \(R\). Observe that \(R \in \Supp(\varepsilon_0^1)(Q)\) if and only if \(\alpha_Q(Y) \in \mpp_R\) and \(R \in \Supp(\varepsilon_1^1)(Q)\) if and only if \(\alpha_Q(X) \in \mpp_R\). Then, the result simply follows from the observation that the ideal \((\alpha_{Q'}(X),\alpha_{Q}(Y))\) contains \((\alpha(X),\alpha(Y))\) and is contained in some maximal ideal of \(F_P[X,Y]\). Such a maximal ideal then corresponds to a place \(R \in M_{F^{\otimes 3}}\) such that \(R_1 = Q\) and \(R_0 = Q'\).
\end{proof}

    For \(n \geq 0\) and \(S \subset \Pl(k)\), we let \(S^{(n)}\) be the subset of \(\Pl(F^{\otimes n+1})\) of places lying above some \(P \in S\). We let \(S_r\) be the set of places of \(k\) that ramify in \(F\). Applying \Cref{prop:TPUnramified} by induction, we get the following:
    \begin{lemma}\label{lemma:UnramifiedPlaces}
        Let \(Q \in \Pl(F^{\otimes n+1}) \setminus S_r^{(n)}\). Then for \(0 \leq i \leq n+1\), we have \(e_{Q,i} = 1\).
    \end{lemma}

    We may now prove the following generalisation of Hilbert's Theorem 90:
\begin{lemma}\label{lemma:Hilbert90}
    Let \(D = \sum_{Q \in \Pl(F^{\otimes 2})} n_Q Q \in \Ker \Dc(\Delta_{Am}^1)\) be supported by places outside of \(S_r^{(1)}\). Then, there exists \(E \in \Dc(F)\) such that \(D = \Dc(\Delta_{Am}^0)(E)\).
\end{lemma}

\begin{proof}
    We set
    \[E = \sum_{P \in \Pl(F)} \left(\min_{\substack{Q \in \Pl(F^{\otimes 2}) \\ Q_0 = P}} n_Q \right) P.\]
    Then, by \Cref{lemma:UnramifiedPlaces}, we get 
    \[\Dc(\epsilon_0^0)(E) = \sum_{Q \in \Pl(F^{\otimes 2})} \left(\min_{\substack{Q'' \in \Pl(F^{\otimes 2}) \\ Q''_0 = Q_0}}  n_{Q''}\right) Q\]
    and
    \[\Dc(\epsilon_1^0)(E) = \sum_{Q \in \Pl(F^{\otimes 2})} \left(\min_{\substack{Q' \in \Pl(F^{\otimes 2}) \\ Q'_0 = Q_1}}  n_{Q'}\right) Q\]

    It follows that
    \[D + \Dc(\epsilon_1^0)(E) = \sum_{Q \in \Pl(F^{\otimes 2})} \left(\min_{\substack{Q' \in \Pl(F^{\otimes 2}) \\ Q'_0 = Q_1}}  n_Q + n_{Q'}\right) Q.\]

    We introduce the following automorphisms:
    \[\begin{array}{rccc}
        \sigma\colon &F^{\otimes 2} &\to &F^{\otimes 2} \\
        & a \otimes b &\mapsto &b \otimes a; \\
        \tau\colon &F^{\otimes 3} &\to &F^{\otimes 3} \\
        &a \otimes b \otimes c &\mapsto &a \otimes c \otimes b.
    \end{array}\]

    Observe that we have the following equalities:
    \begin{align}
        \tau \circ \varepsilon_0^1 &= \varepsilon_0^1 \circ \sigma \label{eq:tau00sigma}\\
        \tau \circ \varepsilon_2^1 &= \varepsilon_1^1 \label{eq:tau21}\\
        \varepsilon_1^0 &= \sigma \circ \varepsilon_0^0 \label{eq:1sigma0}.
    \end{align}
    If \(Q \in \Pl(F^{\otimes 2})\) (resp. \(R \in \Pl(F^{\otimes 3})\)), we write \(Q^\sigma\) (resp. \(R^\tau\)) for the unique place in the support of \(\Dc(\sigma)(Q)\) (resp. \(\Dc(\tau)(R)\)).

    Now, let \(Q,Q' \in \Pl(F^{\otimes 2})\) such that \(Q'_0 = Q_1\). By \Cref{eq:1sigma0}, we have \(Q'_0 = Q_0^\sigma\). Applying \Cref{lemma:CocycleTransit} to \(Q^\sigma\) and \(Q'\), we get \(R \in \Pl(F^{\otimes 3})\) such that \(R_1 = Q'\) and \(R_0 = Q^\sigma\). By \Cref{eq:tau21}, we have \(R^\tau_2 = R_1 = Q'\) and by \Cref{eq:tau00sigma}, we get \(R^\tau_0 = (R_0)^\sigma = Q\). We set \(Q'' = R_1^\tau\). Now, the coefficient of \(R^\tau\) in \(\Dc(\Delta^1_{Am})(D)\) is \(n_Q + n_{Q'} - n_{Q''}\), and since this divisor is zero by hypothesis, we get \(n_{Q''} = n_Q + n_{Q'}\) (note that \(R \notin S_r^{(2)}\), so \(e_{R,i} = 1\) for \(0 \leq i \leq 2\)). Observe that we have
    \begin{align*}
        R &\in \Supp(\Dc(\varepsilon_0^1)(Q)) \\
            &\subset \Supp(\Dc(\varepsilon_0^1 \circ \varepsilon_0^0)(Q_0)) \\
            &= \Supp(\Dc(\varepsilon_1^1 \circ \varepsilon_0^0)(Q_0)) \\
    \end{align*}
    and
    \[R \in \Supp(\Dc(\varepsilon_1^1)(Q'')) \subset \Supp(\Dc(\varepsilon_1^1 \circ \varepsilon_0^0)(Q''_0)).\]
    By \Cref{prop:UniquePlaceAbove}, it follows that \(Q_0 = Q''_0\). That is, we proved that there exists \(Q'' \in \Pl(F^{\otimes 2})\) such that \(Q''_0 = Q_0\) and \(n_Q + n_{Q'} = n_{Q''}\).

    Conversely, let \(Q,Q'' \in \Pl(F^{\otimes 2})\) such that \(Q_0 = Q''_0\). Then by \Cref{lemma:CocycleTransit}, there is \(R \in \Pl(F^{\otimes 3})\) such that \(R_0 = Q\) and \(R_1 = Q''\). Set \(Q' = R_2\) and, as above, we get \(n_Q + n_{Q'} = n_{Q''}\). As above, we use the fact that \(\varepsilon_0^1 \circ \varepsilon_1^0 = \varepsilon_2^1 \circ \varepsilon_0^0\) to prove that \(Q_1 = Q'_0\).

    Putting things together, we have proved that for \(Q \in \Pl(F^{\otimes 2}) \setminus S_r^{(1)}\),
        \[\min_{\substack{Q' \in \Pl(F^{\otimes 2})\\Q'_0 = Q_1}} n_Q + n_{Q'} = \min_{\substack{Q'' \in \Pl(F^{\otimes 2})\\Q''_0 = Q_0}} n_{Q''}.\]

    It follows directly that \(D + \Dc(\varepsilon_1^0)(E) = \Dc(\varepsilon_0^0)(E)\).
    \end{proof}

We now get our main theorem for this section:
\begin{theorem}\label{thm:SUnitTriv}
    Let \(b \in B_{Am}^2(k,F)\) be a coboundary. Let \(S\) be a set of non-archimedean places of \(k\) such that:
    \begin{itemize}
        \item \(S_r \subset S\). That is, \(S\) contains the places of \(k\) that ramify in \(F\).
        \item The finite places of \(S^{(0)}\) generate the class group \(\Cl(F)\).
        \item \(\Supp(\Dc(b)) \subset S^{(2)}\).
    \end{itemize}
    Then there exists a cochain \(a\) in the group of \(S^{(2)}\)-units of \(F^{\otimes 2}\) such that \(b = \Delta_{Am}^1(a)\)
\end{theorem}

\begin{proof}
 Let \(\alpha \in (F^{\otimes 2})^\times\) be such that \(\Delta_{Am}^1(\alpha) = b\). We consider the divisor \(D = \Dc(\alpha) = \sum_{Q \in \Pl(F^{\otimes 2})} n_Q Q\) of \(\alpha\). We set \(D_S = \sum_{Q \in S^{(1)}} n_Q Q\) and \(D'_{S} = \sum_{Q \notin S^{(1)}} n_Q Q\). Now we get \(\Dc(\Delta_{Am}^1)(D) = \Dc(b)\) and therefore it is supported by \(S^{(2)}\). Observe that if \(Q \in S^{(1)}\), then \(\Supp(\Dc(\Delta_{Am}^1)(Q)) \subset S^{(2)}\) and conversely, if \(\Supp(\Dc(\Delta_{Am}^1)(Q)) \subset S^{(2)}\) and \(\Dc(\Delta_{Am}^1)^1(Q) \neq 0\), then \(Q \in S^{(1)}\). It follows that \(\Dc(\Delta_{Am}^1)(D'_S) = \Dc(\Delta_{Am}^1)(D) - \Dc(\Delta_{Am}^1)(D_S) = 0\).

 The support of \(D'_{S}\) is disjoint from \(S_r^{(1)}\). We may therefore apply \Cref{lemma:Hilbert90} and get a divisor \(E \in \Dc(F)\) such that \(D'_S = \Delta_{Am}^0(E)\). Now, as \(S^{(0)}\) generates the class group of \(F\), there exists \(E' \in \Dc(F)\) with support in \(S^{(0)}\) and \(\gamma \in F^\times\) such that \(E = \Dc(\gamma) + E'\). Then, we get that
 \[\Dc(\Delta_{Am}^0)(\Dc(\gamma)) + \Dc(\Delta_{Am}^0)(E') = D - D_S\]
 and therefore
 \[\Dc(\Delta_{Am}^0)(E') + D_{S} = \Dc(\alpha \Delta_{Am}^0(\gamma^{-1})).\]
 Now, since \(\Dc(\Delta_{Am}^0)(E') + D_S\) has support in \(S^{(1)}\), this shows that \(\alpha \Delta_{Am}^0(\gamma^{-1})\) is a \(S^{(1)}\)-unit. Furthermore, \[\Delta_{Am}^1(\alpha\Delta_{Am}^0(\gamma^{-1})) = \Delta_{Am}^1(\alpha) = b,\]
 and \(\alpha \Delta_{Am}^0(\gamma^{-1})\) is a cochain with the required properties.
\end{proof}

From \Cref{thm:SUnitTriv} we directly get an algorithm for computing a trivialisation of a \(2\)-coboundary:

\begin{algorithm}
    \KwIn{A number field \(k\)}
    \KwIn{An étale \(k\)-algebra \(F\) defined by a separable polynomial \(P \in k[X]\)}
    \KwIn{A coboundary \(b \in B_{Am}^2(k,F)\)}
    \KwOut{A cochain \(a \in C_{Am}^0(k,F)\) such that \(\Delta_{Am}^1(a) = c\)}
    Compute \(S_1\), the set of places of \(k\) that ramify in \(F\)\;
    Compute \(S_2\), a set of places of \(k\) such that \(S_2^{(0)}\) generates the class group \(\Cl(F)\)\;
    Compute \(\Dc(b)\). Let \(S_3\) be the set of places of \(k\) below the places in the support of \(\Dc(b)\)\;
    Set \(S = S_1 \cup S_2 \cup S_3\)\;
    Compute the sets \(S^{(2)}\) and \(S^{(3)}\)\;
    Compute an isomorphism \(\phi\) from the group of \(S^{(2)}\)-units of \(F^{\otimes 2}\) to \(\Z^r \oplus \Z/m\Z\)\;
    Compute an isomorphism \(\psi\) from the group of \(S^{(3)}\)-units of \(F^{\otimes 3}\) to \(\Z^{r'} \oplus \Z/m'\Z\)\;
    Solve the linear equation \((\psi \circ \Delta_{Am}^1 \circ \phi^{-1})(a) = \psi(b)\)\;
    \KwRet{\(a\)}
    \caption{Computing a trivialisation of a \(2\)-coboundary}
    \label{algo:TrivCobound}
\end{algorithm}

\begin{theorem}\label{thm:AlgoTrivCobound}
    Given a number field \(k\), a separable polynomial \(P \in k[X]\)  and a coboundary \(b \in B^2_{Am}(k,F)\), where \(F = k[X]/P\), \Cref{algo:TrivCobound} outputs the representation of a cochain \(\alpha \in C^1(k,F)\) such that \(\Delta_{Am}^1(\alpha) = b\). Furthermore, assuming GRH, \Cref{algo:TrivCobound} is a polynomial quantum algorithm.
\end{theorem}

\begin{proof}
    Using a polynomial-time algorithm for factoring polynomials over number fields \cite{lenstra1983factoring}, one may compute splittings of \(F\), \(F^{\otimes 2}\) and \(F^{\otimes 3}\) as direct products of number fields. Set \(F^{\otimes n+1} = F^{(n)}_1 \times \hdots \times F_{r_n}^{(n)}\).

    The set \(S_1\) may be computed by factoring the discriminants of the extensions \(F^{(0)}_i/k\), which may be done in polynomial time by \Cref{fact:maxord}. Then, by \Cref{fact:Bach}, one may compute in polynomial time a set of polynomial size of places of \(F\) that generate the class group \(\Cl(F)\), and the set \(S_2\) of places of \(k\) lying below them. In practice, one may set \(S_2\) to be the set of prime ideals \(\p\) of \(k\) such that \(N(\p) \leq 12\log^2(\max_{i \in [r]} |\Delta_{F_i}|)\). We may then take the image \((b_1,\hdots,b_{r_2})\) of \(b\) in the product \(F^{(2)}_1 \times \hdots \times F^{(2)}_{r_2}\) and factor the principal ideals \((b_i)\) in \(F^{(2)}_i\) into a product \(\p_{i,1}\hdots \p_{i,s_i}\), which may be done in polynomial time by \Cref{fact:factoring}. We may then set \(S_3 = \{\p_{i,j} \cap k, i \in [r_2], j \in [s_i]\} \subset \Pl(k)\).
    
   Isomorphisms \(\phi\) and \(\psi\) are computed using \Cref{fact:s-unit} and \cite[Algorithm 8.4]{lenstra2018algorithms}. Finally, the last step is the computation of a solution of a system of linear equations over \(\Z\).
   
    The correctness of the algorithm relies on the fact that a cochain \(a\) such that \(b = \Delta_{Am}^1(a)\) exists and may be found in the group of \(S^{(1)}\)-units of \(F^{\otimes 2}\), which is the content of \Cref{thm:SUnitTriv}.
\end{proof}

\begin{theorem}\label{cor:QuantumSplit}
    Assuming GRH, there exists a polynomial quantum algorithm which, give a number field \(k\) and an algebra \(A \simeq M_n(k)\), computes an explicit algorithm from \(A\) to \(M_n(k)\).
\end{theorem}

\begin{proof}
    This is simply a combination of \Cref{thm:AlgoFindCocycle,thm:AlgoTrivCobound}. Indeed, using \Cref{algo:FindCocycle}, one may compute an étale \(k\) algebra \(F = k[X]/P\), a cocycle \(c \in \Z^2(k,F)\) and an isomorphism \(A \simeq A(F,c)\). Since \(A\) is isomorphic to \(M_n(k)\), the cocycle \(c\) is in fact a coboundary. Then, a trivialisation of \(c\) may be computed using \Cref{algo:TrivCobound}. Applying \Cref{cor:AssoCocyIsomAmitsur}, we obtain an explicit isomorphism \(A(F,c) \simeq A(F,\1)\). Finally, an isomorphism \(A(F,\1) \simeq M_n(k)\) may easily be computed using \Cref{remark:AmitsurTrivialCocycleIsom}.
\end{proof}




\bibliographystyle{acm}
\bibliography{biblio}
\end{document}